\documentclass[12pt]{article}
\usepackage{hyperref}
\usepackage{amsfonts}
\usepackage{enumerate}
\begin{document}

\begin{center}
{\large\bf \uppercase{Controllability and reachability of singular linear discrete time systems}}

\vskip.20in
Charalambos P. Kontzalis$^{1}$ and\ Grigoris Kalogeropoulos$^{2}$\\[2mm]
{\footnotesize
$^{1}$Department of Informatics, Ionian University, Corfu, Greece\\[5pt]
$^{2}$Department of Mathematics, University of Athens, Greece}
\end{center}

{\footnotesize
\noindent
\textbf{Abstract:} The main objective of this article is to develop a matrix pencil approach for the study of the controllability and reachability of a class of linear singular discrete time systems. The description equation of a practical system may be established through selection of the proper state variables. Time domain analysis is the method of analyzing the system based on this description equation, through which we may gain a fair understanding of the system's structural features as well as its internal properties. Using time domain analysis, this article studies the fundamentals in system theory such as reachability and controllability. \\
\\[3pt]
\\[3pt]

\vskip.2in
\section{Introduction}

Linear discrete time systems are systems in which the variables take their value at instantaneous time points. Discrete time systems differ from continuous time ones in that their signals are in the form of sampled data. With the development of the digital computer, the discrete time system theory plays an important role in control theory. In real systems, the discrete time system often appears when it is the result of sampling the continuous-time system or when only discrete data are available for use. Discrete time systems have many applications in economics, physics, circuit theory, and other areas. For example in finance, there is the very famous Leondief model, see [2], or the very important Leslie population growth model and backward population projection, see also [2]. In physics the Host-parasitoid Models, see [46]. Applications of absorbing Markov chains or the distribution of heat through a long rod or bar are other interesting applications suggested in [46]. Thus many authors have studied discrete time systems, see and their applications, see [1-12, 15-21, 24-34, 39-46, 48-50]. In this article we study the controllability and reachability of these systems and extend known results in the literature. The results of this paper can be applied also in systems of fractional nabla difference equations, see [13, 14, 22, 35-38]. In addition they are very useful for applications in many mathematical models using systems of difference equations existing in the literature, see [1, 2, 16-18, 46, 48, 49]. We consider
\begin{equation}
FX_{k+1}=GX_k+BV_k
\end{equation}
where  $F \in \mathcal{M}({m \times m;\mathcal{F}})$ is a singular matrix (det$F$=0), $G, B\in \mathcal{M}({m \times m;\mathcal{F}})$ and $X_k, V_k\in \mathcal{M}({m \times
1;\mathcal{F}})$. $\mathcal{M}({m \times m;\mathcal{F}})$ is the
algebra of square matrices with elements in the field
$\mathcal{F}$). For the sake of simplicity we set
$\mathcal{M}_m  = \mathcal{M}({m \times m;\mathcal{F}})$, $\mathcal{M}_{nm}  = \mathcal{M}({m \times n;\mathcal{F}} )$. In the following sections we will study the controllability and reachability of these systems.

\section{Mathematical background and notation}

This brief subsection introduces some preliminary concepts and
definitions from matrix pencil theory, which are being used
throughout the paper. Linear systems of type (1) are closely
related to matrix pencil theory, since the algebraic geometric, and
dynamic properties stem from the structure by the associated pencil
$sF-G$. The matrix pencil theory has been extensively used for the study of LDTSs with time invariant
coefficients, see for instance [7-21, 23-27, 32-34, 43]. The class of $sF-G$ is characterized by a uniquely
defined element, known as a complex Weierstrass canonical form,
$sF_w -G_w$, see [7, 23, 32, 43], specified by the complete set of
invariants of $sF-G$. This is the set of \emph{elementary divisors} (e.d.) of the following type:
\begin{itemize}
    \item e.d. of the type  $(s-a_i)^p_i$, \emph{are called finite elementary
    divisors} (nz. f.e.d.) and $a_i$ are the finite eigenvalues with multiplicity $p_i$.
    \item e.d. of the type  $\hat{s}^q$ are called \emph{infinite elementary divisors}
    (i.e.d.), of multiplicity $q$.
\end{itemize}
Then, the complex Weierstrass form $sF_w -G_w$  of the regular
pencil $sF-G$ is defined by 
\[
sF_w  - Q_w = blockdiag (sI_p  - J_p, sH_q  - I_q),
\]
where $p+q=m$ and the first normal Jordan type element is uniquely defined by the set of f.e.d. (the $p$ finite eigenvalues),
\[
  ({s - a_1 })^{p_1 } , \dots ,({s - a_\nu  }
 )^{p_\nu },\quad \sum_{j = 1}^\nu  {p_j  = p}
\]
of $sF-G$  and has the form
\[
    sI_p  - J_p  = blockdiag (sI_{p_1 }  - J_{p_1 } (
    {a_1 }), \dots, sI_{p_\nu  }  - J_{p_\nu  }
    ({a_\nu  }) )
\]
And also the $q$  blocks of the second uniquely defined block
$sH_q -I_q$ correspond to the i.e.d.
\[
  \hat s^{q_1} , \dots ,\hat s^{q_\sigma}, \quad \sum_{j =
  1}^\sigma  {q_j  = q}
\]
of $sF-G$ and has the form
\[
    sH_q  - I_q  = blockdiag (sH_{q_1 }  - I_{q_1 },\dots , sH_{q_\sigma  }  - I_{q_\sigma})
\]
$J_{p_j } ({a_j }),H_{q_j }$ are defined as
\[
   J_{p_j } ({a_j }) =  \left[\begin{array}{ccccc}
   a_j  & 1 & \dots&0  & 0  \\
   0 & a_j  &   \dots&0  & 0  \\
    \vdots  &  \vdots  &  \ddots  &  \vdots  &  \vdots   \\
   0 & 0 &  \ldots& a_j& 1\\
   0 & 0 & \ldots& 0& a_j
   \end{array}\right] \in {\mathcal{M}}_{p_j }
\]
\[
 H_{q_j }  = \left[
\begin{array}{ccccc} 0&1&\ldots&0&0\\0&0&\ldots&0&0\\\vdots&\vdots&\ddots&\vdots&\vdots\\0&0&\ldots&0&1\\0&0&\ldots&0&0
\end{array}
\right] \in {\mathcal{M}}_{q_j }.
  \]
There exist nonsingular $\mathcal{M}(m\times m,\mathcal{F})$ matrices $P$ and $Q$ such that 
\[
\begin{array}{c}PFQ = F_w  = blockdiag(I_p, H_q)\\
PGQ = G_w  = blockdiag(J_p, I_q).
\end{array}
\]
Let 
\[
PB=\left[\begin{array}{c}
 B_p\\ B_q
 \end{array}\right]
\] 
where $B_p\in\mathcal{M}_{pm}, B_q\in\mathcal{M}_{qm}$
\\\\
\textbf{Lemma 2.1.}
System (1) is divided into two subsystems, the subsystem
\begin{equation}
    Y^p_{k+1} =  J_p Y^p_k+B_p V_k, 
\end{equation}
and the subsystem
\begin{equation}
    H_q Y^q_{k+1} = Y^q_k+B_q V_k
\end{equation}
\textbf{Proof.}
Consider the transformation
\[
    X_k=QY_k
\]
Substituting the previous expression into (1) we obtain
\[
    FQY_{k+1}=GQY_k+BV_k
\]
Whereby, multiplying by P, we arrive at
\[
    F_wY_{k+1}=G_w Y_k+PBV_k
\]
Moreover, we can write 
\[
Y_k =\left[\begin{array}{c}
 Y^p_k \\
 Y^q_k
 \end{array}\right]
\]
taking into account the above expressions, we arrive easily at
(2) and (3).
\\\\
For the sake of simplicity, system (1) is assumed to be in its standard decomposition from (2) and (3). However the results are applicable to systems in the general form of (1). We start from the concept of reachable set to study the state structure. Then from [7-21, 23-27, 32-34, 43] we know that the state response for systems (2), (3) are 
\[
    Y^p_k=J_p^kY^p_0+\sum^{k-1}_{i=0}J_p^{k-i-1}B_pV_i , k\geq 0, 
\]
and
\[
Y^q_k=-\sum^{q_*-1}_{i=0}H_q^iB_qV_{k+i}
\]
The reachable set may be defined in the next chapter.

\section{Reachability}

\textbf{Definition 3.1.} Any vector $W\in\mathcal{M}_{m1}$ in m dimensional vector space is said to be reachable from an initial condition $Y_0$, if there exists a sequence of inputs $V_0, V_1, ..., V_k, ...$, k=0, 1, ..., such that $Y_k=W$
\\\\
For a given system we introduce the following notation
\\\\
$\Re Y_0$ =(the set of all states reachable from $Y_0$)=Y$\in\mathcal{M}_{m1}$: there exists a sequence of inputs $V_0, V_1, ..., V_k, ...$ such that $Y_k=Y$ which is a subspace in a m vector space.
\\\\
If for i, j= 1, 2, ... we define
\[
S_i(J_p/B_p)=colspan\left[\begin{array}{cccc}B_p&J_pB_p&\cdots&J^{i-1}_pB_p\end{array}\right]
\]
and
\[
S_j(H_q/B_q)=colspan\left[\begin{array}{cccc}B_q&H_qB_q&\cdots&H^{j-1}_qB_q\end{array}\right]
\]
Also we define 
\[
\bar S_p(J_p/B_p)=\left\{\begin{array}{cccc}Y\in M_{m1}:&Y=\left[\begin{array}{c}Y^p_1\\Y^q_2\end{array}\right]&Y^p_1\in S_i((J_p/B_p)&Y^q_2 =\bar 0 \end{array}\right\}
\]
and
\[
\bar S_q(H_q/B_q)=\left\{\begin{array}{cccc}Y\in M_{m1}:&Y=\left[\begin{array}{c}Y^p_1\\Y^q_2\end{array}\right]&Y^p_1=\bar 0&Y^q_2 \in S_q((H_q/B_q) \end{array}\right\}
\]
Then we can state the following theorem.
\\\\
\textbf{Theorem 3.1.} The state reachable set $\Re(\bar0)$ from the zero initial condition $Y_0=\bar 0$ is given by:
\[
\Re (\bar 0)=\bar S_p(J_p/B_p)\oplus\bar S_q^*(H_q/B_q)
\]
where $q^*$ is the index of annihilation of $H_q$.
\\\\
\textbf{Proof.} Solving the subsystems (2), (3) we have respectively, see [7-21, 23-27, 32-34, 43], the following solutions:
\[
 Y^p_k=J_p^kY^p_0+\sum^{k-1}_{i=0}J_p^{k-i-1}P_1V_i 
\]
and
\[
  Y^q_k=-\sum^{q_*-1}_{i=0}H_q^iP_2V_{k+i}.
\]
Setting $Y_0=\bar 0$ in (2) we have:
\[
Y^p_k=\sum^{k-1}_{j=0}J^{k-1-j}_pB_pV_j=J^{k-1}_pB_pV_j+J^{k-2}_pB_pV_j+...+B_pV_{k-1},
\]
or in matrix form
\[
Y^p_k=\left[\begin{array}{cccc}B_p&J_pB_p&\cdots&J^{i-1}_pB_p\end{array}\right]\left[\begin{array}{c}V_{k-1}\\V_{k-2}\\...\\V_{k_0}\end{array}\right].
\]
We consider the following two cases:
\\\\
(a) $k\leq p$, then it is obvious that for every $Y_k \in S_k(J_p/B_p)$ we have 
\[
Y_k \in S_p(J_p/B_p)
\]
(b) $k>p$, then if $Y_k\in S_k(J_p/B_p)$ we have, 
\[
Y_k=\sum^{p}_{i=1}J^{i-1}_pB_pC_i+\sum^{k}_{j=p+1}J^{j-1}_pB_pC_j
\]
For $j\geq p+1$ every power $J_p^{j-1}$, can be represented as a linear sum of powers of $J^0_p, J^1_p, ..., J^{p-1}_p$ (this is a straight consequence of the Cayley-Hamilton theorem) and obviously 
\[
Y_k\in S_p(J_p/B_p)
\]
So for every $k\geq 0$ we have 
\[
S_k(J_p/B_p)\subseteq S_p(J_p/B_p)
\]
Let's now take system (3), then expanding the sum giving $Y^q_k$ we get:
\[
Y^q_k=-B_qV_k-H_qB_qV_{k+1}-...-H^{q^*-1}B_qV_{k+q*-1}
\]
or in matrix form
\[
Y^q_k=\left[\begin{array}{cccc}B_q&H_qB_q&\cdots&H^{q*-1}_qB_q\end{array}\right]\left[\begin{array}{c}-V_k\\-V_{k+1}\\...\\-V_{k+q*-1}\end{array}\right]\in S_{q^*}(H_q/B_q)
\]
Then every state $Y_k$, k starting from zero, is of the form:
\[
Y_k=\left[\begin{array}{c}Y^p_k\\Y^q_k\end{array}\right]=\left[\begin{array}{c}Y^p_k\\\bar 0\end{array}\right]+\left[\begin{array}{c}\bar0\\Y^q_k\end{array}\right]\in S_p(J_p/B_p)\oplus S_q^{q*}(H_p/B_q)
\]
The last observation leads to the fact that 
\begin{equation}
\Re (\bar0)\subseteq S_p(J_p/B_p)\oplus S_q^{q*}(H_p/B_q)
\end{equation}
where by $\Re (\bar0)$ we denote the reachable set from zero of the system (1).
\\\\
On the other hand, if we take $Y^*=\left[\begin{array}{c}Y^*_p\\Y^*_q\end{array}\right]\in \bar S_p(J_p/B_p)\oplus \bar S_q^{q*}(H_p/B_q)$ we will prove that there exists a sequence of inputs $V_0, V_1, V_2, ...$ such that the trajectory starting from zero  $Y_0=\bar 0$, arrives in $Y^*$ in k-steps. In this case $Y^*_p\in S_p(J_p/B_p)$ and $Y^*_q \in S_q^{q*}(H_p/B_q)$. We will see that there exists a sequence of inputs $V_0, V_1, V_2, ...$ and a suitable k such that:
\[
Y^*_p=\sum^{k-1}_{j=0}J^{k-1-j}_pB_pV_j
\]
and
\[
Y^*_q=-\sum^{q_*-1}_{i=0}H_q^iB_qV_{k+i}
\]
It is known that, for any $q^*$ there always exists a vector polynomial $f_k$, whose order is $q^*-1$, such that
\[
f_{k+i}=V_i,i=0, 1, 2, ..., q^*-1,
\]
for the proof see the Higher Degree Interpolation theorem [7]. If we impose the input control $V_k=V^1_k+f_k$, then for $Y=\left[\begin{array}{c}Y^p_1\\Y^q_2\end{array}\right]\in\Re(\bar 0)$ it will be that 
\begin{equation}
Y^p_1=\sum^{k-1}_{j=0}J^{k-1-j}_pB_pV^1_k+\sum^{k-1}_{j=0}J^{k-1-j}_pB_pf_k
\end{equation}
Let
\[
\bar Y^p_1=Y^*_p-\sum^{k-1}_{j=0}J^{k-1-j}_pB_pf_k.
\]
In this case there exists input sequence $V^1_k$, see  [7], such that 
\[
\bar Y^p_1=\sum^{k-1}_{j=0}J^{k-1-j}_pB_pV^1_k
\]
and
\[
\sum^{k-1}_{i=0}H^i_qB_qV_{k+i}=\bar 0
\]
From (5) we have that 
\[
Y^p_1=\bar Y^p_1+Y^*_p-\bar Y^p_1=Y^*_p
\]
and
\[
Y^q_2=-\sum^{q*-1}_{i=0}H^i_qB_qV^1_{k+i}-\sum^{q*-1}_{i=0}H^{i}_qB_qf_{k+i}=\bar 0-\sum^{q*-1}_{i=0}H^{i}_qB_qV_i=Y^*_q.
\]
Thus
\begin{equation}
\Re (\bar0) \supseteq S_p(J_p/B_p)\oplus S_q^{q*}(H_p/B_q)
\end{equation}
and from (4), (6) the proof is complete.
\section{Finite discrete time systems}
Consider the finite discrete time system described by
\begin{equation}
\left\{\begin{array}{cc}FX_{k+1}=GX_k+BV_k,&0\leq k\leq M\\X_k=QY_k\end{array}\right\}
\end{equation}
where $F, G, B, Q \in\mathcal{M}_m$, $X_k, Y_k \in \mathcal{M}_{m1}$, F is singular and the system pencil sF-G is regular. From Lemma 2.1, the system (7) can be divided into the subsystems (2), (3). From [7-21, 23-27, 32-34, 43], the solution of system (2) is given by 
\begin{equation}
    Y^p_k=J_p^kY^p_0+\sum^{k-1}_{i=0}J_p^{k-i-1}B_pV_i , k\geq 0.
\end{equation}
Using the terminal point $Y^q_M$ the solution of the system (3) is given by the following proposition.
\\\\
\textbf{Proposition 4.1.} Consider the system (3). Then for $0\leq k\leq M$, the solution is given from
\begin{equation}
Y^q_k=H^{M-k}_qY^q_M-\sum^{M-k-1}_{i=0}H^i_qB_qV_{k+i}.
\end{equation}
\textbf{Proof.} From the equation (3) we can obtain the following equations
\[
\begin{array}{c}
    H_qY_{k+1}^q=Y_k^q+B_qV_k\\
    H_q^2Y_{k+1}^q=H_qY_k^q+H_qB_qV_k\\
    H_q^3Y_{k+1}^q=H_q^2Y_k^q+H_q^2B_qV_k\\
    \vdots\\
    H_q^{M-q}Y_{k+1}^q=H_q^{m-k-1}Y_k^q+H_q^{M-k-1}B_qV_{M-1}
    \end{array}
    \]
    and
    \[
		\begin{array}{c}
    H_qY_{k+1}^q=Y_k^q+B_qV_k\\
    H_q^2Y_{k+2}^q=H_qY_{k+1}^q+H_qB_qV_{k+1}\\
    H_q^3Y_{k+3}^q=H_q^2Y_{k+2}^q+H_q^2B_qV_{k+2}\\
    \vdots\\
    H_q^{M-q}Y_M^q=H_q^{M-k-1}Y_{M-1}^q+H_q^{M-k-1}B_qV_{M-1}.
    \end{array}
    \]
by repetitively substitution of each equation in the next one, the conclusion
\[
Y^q_k=H_q^{M-k}Y^q_M-\sum^{M-k-1}_{i=0}H^i_qB_qV_{k+i}
\]
is obtained.
\\\\
\textbf{Remark 4.1.} As we have seen above, the state at any time point k for a finite discrete time system is generally related with not only the initial state and former inputs, as in the regular system case, but also terminal state $X_M$ an future inputs up to time point $M$.
\\\\
\textbf{Remark 4.2.} It is clear that the solution of system (1) is depending on the initial point $Y^p_0$ and former inputs $V_0, V_1, ..., V_{k-1}$ but the solution of the system (7), because of the solution (9) of subsystem (3), is depending on the terminal point $Y_M^q$ and future inputs.
\\\\
\textbf{Remark 4.3.} The solution of the regular system (2) is given by the formula (8). So the state $Y_k^p$ is uniquely determined by the initial state $Y_0^p$ and $V_k$, with $0\leq k\leq M$.
\\\\
\textbf{Example 4.1.} Consider the finite discrete time system for $0\leq k\leq M$ in canonical form given by
\[
\left[\begin{array}{cccc}1&0&0&0\\0&1&0&0\\0&0&1&0\\0&0&0&0\end{array}\right]Y_{k+1}=\left[\begin{array}{cccc}1&1&0&0\\0&1&0&0\\0&0&1&0\\0&0&0&1\end{array}\right]Y_k+\left[\begin{array}{c}0\\1\\1\\-1\end{array}\right]V_k
\]
By writing $Y_k=\left[\begin{array}{c} Y^p_k\\Y^q_k \end{array} \right]$ the system maybe written in the form of two subsystems 
\[
Y^p_{k+1}=\left[\begin{array}{cc}1&1\\0&1\end{array}\right]Y^p_k+\left[\begin{array}{c}0\\1\end{array}\right]V_k
\]
and
\[
\left[\begin{array}{cc}0&1\\0&0\end{array}\right]Y^q_{k+1}=Y^q_k+\left[\begin{array}{cccc}1\\-1\end{array}\right]V_k
\]
From (14) and (31) its solution for $0\leq k\leq M$ is given by
\[
Y^p_k=\left[\begin{array}{cc}1&1\\0&1\end{array}\right]^kY^p_0+\sum^{k-1}_{i=0}\left[\begin{array}{c}k-i-1\\1\end{array}\right]V_i
\]  
or
\[
Y^p_k=\left[\begin{array}{cc}1&k\\0&1\end{array}\right]Y^p_0+\sum^{k-1}_{i=0}\left[\begin{array}{c}k-i-1\\1\end{array}\right]V_i
\] 
and 
\[
Y^q_k=\left[\begin{array}{cc}0&1\\0&0\end{array}\right]^{M-k}Y^q_M-\sum^{M-k-1}_{i=0}\left[\begin{array}{cc}0&1\\0&0\end{array}\right]^i\left[\begin{array}{c}-1\\1\end{array}\right]V_{k+i}
\]
or
\[
Y^q_k=\begin{array}{cc}\left[\begin{array}{cc}0&1\\0&0\end{array}\right]Y^q_M-\left[\begin{array}{c}1\\-1\end{array}\right]V_k,&k=M-1\end{array}
\]
and
\[
Y^q_k=\begin{array}{cc}\left[\begin{array}{c}1\\0\end{array}\right]V_{k+1}-\left[\begin{array}{c}1\\-1\end{array}\right]V_k,&0\leq k\leq M-2\end{array}.
\]
Clearly $Y_k^q$ is independent of the terminal state when $k\leq M-2$. In the next section we will introduce the concept of controllability for finite discrete time systems of the form (7).

\section{Controllability}
\textbf{Definition 5.1.} Finite discrete time system of the form (7) is called controllable if and only if for $W\in M_{m1}$ and any complete set of conditions $Y^p_0$ and $Y^q_M$ there exists a time point $k_1$, $0 \leq k_1 \leq M$ and a sequence of control inputs $V_0, V_1, ..., V_M$ such that $Y_{k_1}=W$.
\\\\
Next we will give necessary and sufficient condition for the system of the form (7) to be controllable.
\\\\
\textbf{Theorem 5.1.} The system (7) is controllable if and only if
\[
rank\left[\begin{array}{cccc}B_p&J_pB_p&\cdots&J^{i-1}_pB_p\end{array}\right]=p
\]
and
\[
rank\left[\begin{array}{cccc}B_q&H_qB_q&\cdots&H^{j-1}_qB_q\end{array}\right]=q.
\]
Where $J_p, B_p, B_q, H_q$ are determined in section 2.
\\\\
\textbf{Proof.} \textit{Necessity:} Let $Y^p_0=Y^q_M=\bar 0$. Under the controllability assumption, for any $W\in M_{m1}$ there exists a time point $k_1$ and a sequence input $V_i$, i=0, 1, 2,..., M, such that $Y_{k_1}=W$. On the other hand, from the state representation, we have
\begin{equation}
W=Y_{k_1}=TV
\end{equation}
Where
\[
Y_{k_1}=\left[\begin{array}{c} Y^p_{k_1}\\Y^q_{k_1} \end{array} \right]
\]
\[
T=diag(\left[\begin{array}{cccc}B_p&J_pB_p&\cdots&J^{i-1}_pB_p\end{array}\right], \left[\begin{array}{cccc}B_q&H_qB_q&\cdots&H^{j-1}_qB_q\end{array}\right])
\]
\[
V=\left[\begin{array}{c}V_0\\V_1\\\vdots\\ V_M\end{array}\right]
\]
From (10) and the arbitrariness of W we know that the conclusion holds.\\
\textit{Sufficiency}: For any condition $Y^p_0$, $Y^q_M$ the state $Y_k$ has the form
\begin{equation}
Y_k=TV+U_k
\end{equation}
Where
\[
U_k=\left[\begin{array}{c}J^k_pY^p_0\\H_q^{M-k}Y^q_M\end{array}\right] 
\]
Under the sufficient assumption, the matrix $T$ is full row rank ($M\geq m$). Therefore for any $W$ we choose $k_1=p$ and
\[
\left[\begin{array}{c}V_0\\V_1\\\vdots\\ V_{M-1}\end{array}\right]=T^t(TT^t)^{-1}(W-U_p)
\]
where $^t$ is the transpose tensor. The inputs determined here will satisfy $Y^p_k=W$. Therefore (7) is controllable if and only if its subsystems (2), (3), for $k=0,1,...,M$, are controllable. Since $Y^p_k$ is governed by $V_0, V_1, ..., V_{k-1}$ and $Y^q_k$ by $V_k, V_{k+1}, ..., V_{M-1}$, we may choose control inputs respectively for the control purpose of $Y^p_k$ and $Y^q_k$. From a given fixed terminal condition $Y^q_M$ we use $\Re(Y^q_M)$ to denote the reachable state set of the system (7) and this is defined by
\\\\ 
$\Re(Y^q_M)=[ W \in \mathcal{M}_{m1}$: $\exists Y^p_0$, $0\leq k_1 \leq M$ and $V_0, ..., V_M$ such that $Y_{k_1}=W ]$
\\\\
It is clear that initial reachable set $\Re(Y^q_M)$ is dependent of $Y^q_M$. For different $Y^q_M$, $\Re(Y^q_M)$ may be different.
\\\\
\textbf{Proposition 5.1.} For any fixed terminal condition $Y^q_M$ the initial reachable set $\Re(Y^q_M)$ is given by
\begin{equation}
\Re(Y^q_M)=M_{p1}\oplus[\left\{H_q^{M-k}Y^q_M-\sum^{M-k-1}_{i=0}H^i_qB_qV_{k+i}\right\}\cup\left\{Y^q_M\right\}]
\end{equation}
\textbf{Proof.} The result is a simple consequence of the formulas in (8) and (9) giving solutions and solutions space of the subsystems (2), (3) for $k=0,1,...,M$.
\\\\
\textbf{Example 5.1.} Consider the finite time discrete time system for $0\leq k \leq M$ in canonical form given by
\[
\left[\begin{array}{ccc}1&0&0\\0&1&0\\0&0&0\end{array}\right]Y_{k+1}=\left[\begin{array}{ccc}2&0&0\\0&1&0\\0&0&1\end{array}\right]Y_k+\left[\begin{array}{c}1\\1\\0\end{array}\right]V_k
\]
For any complete condition $\left[\begin{array}{c} Y^p_0\\Y^q_M \end{array} \right]$ its state $Y_k$=$\left[\begin{array}{c} Y^p_k\\Y^q_k \end{array} \right]$ is given by
\[
Y^p_k=2^kY^p_0+\sum^{k-1}_{i=0}2^{k-i-1}V_i
\] 
and 
\[
Y^q_k=\begin{array}{cc}\left[\begin{array}{cc}0&1\\0&0\end{array}\right]Y^q_M,&k=M-1\\0,&0\leq k\leq M-2\end{array}
\]

Following the formula (10) the initial reachable set for any terminal condition $Y^q_M$ is
\[
\Re(Y^q_M)=M_{p1}\oplus[\left\{\left[\begin{array}{cc}0&1\\0&0\end{array}\right]Y^q_M\right\}\cup\left\{Y^q_M\right\}]
\]

\section*{Conclusions}

In this article, we give first the definition of the reachable set from an initial condition for systems of the form of (1) and we compute the form of the reachable set from zero initial condition ($Y_0=\bar 0$). Next we consider the finite discrete time system and we give the solution in explicit form. We observe that the state at any time point $k$ for a finite discrete time system is related not only the initial state and former inputs, but also terminal state and future inputs up to the point M. The definition of controllability for finite discrete time systems is given as well  as the necessary and sufficient conditions for such a system to be controllable. Finally for finite discrete time systems we define the initial reachable set from a fixed terminal condition and we give the description of this set. An example given at the end of the section makes more clear the notion and the computation of the initial reachable set from a given fixed terminal condition. As a further extension of this article is to to study controllability, reachability of systems of fractional nabla difference equations. For all this there is already some research in progress.

\end{document}